\documentclass[12pt]{article}
\usepackage{ifpdf}
\usepackage {graphics}
\usepackage{amsmath,mathptmx,amssymb,bm}%,amssymb,amsthm,enumerate,cite,mathrsfs%,backref

\begin{document}

%\begin{frontmatter}

\title{Exact solutions via equivalence transformations of  variable-coefficient fifth-order KdV equations}

%\author{ \\ Departamento de Matem\'aticas, Universidad de C\'adiz, 11510,
% \\ Puerto Real, C\'adiz, Spain\\\date{ {e.mail: 
% }}}
%\date{\empty}

\author{M.S. Bruz\'on ${}^{a}$, R. de la Rosa ${}^{b}$, R. Tracin\`{a} ${}^{c}$\\
  ${}^{a}$ Universidad de C\'adiz, Spain  (e-mail:  m.bruzon@uca.es). \\
 ${}^{b}$ Universidad de C\'adiz, Spain  (e-mail: rafael.delarosa@uca.es). \\
  ${}^{c}$ Universit\`a di Catania, Italy (e-mail:  tracina@dmi.unict.it). \\
}

\date{}
 
\maketitle

\begin{abstract}

In this paper, a family of variable-coefficient fifth-order KdV equations has been considered. By using an infinitesimal method based on the determination of the equivalence group, differential invariants and invariant equations are obtained. Invariants provide an alternative way to find equations from the family which may be equivalent to  a specific subclass of the same family and the invertible transformation which maps both equivalent equations. Here, differential invariants are applied to obtain exact solutions.\\


%\noindent \textit{PACS}: 02.20.Hj; 02.20.Sv; 02.30.Em; 02.30.Jr; 47.10.ab.\\
%\noindent \textit{MSC}: 35C07; 35Q53.
\end{abstract}

%\MOS{XXX; XXX}

%\keywords{Lie Symmetries; Nonlinear self-adjointness; Conservation
%laws; Multipliers}

\maketitle

\section{Introduction}\label{intro}

In this paper, we are interested to look for exact solutions  via equivalence transformations of fifth-order KdV equations with time-dependent coefficients and linear
damping term of the form
\begin{equation}\label{ed1}
\begin{array}{l}
 u_t+A(t) u_{xxxxx}+B(t) u_{xxx}+ C(t) u u_{xxx} +E(t) u u_x+ F(t) u_x u_{xx} + Q(t) u=0 , 
\end{array}
\end{equation}
where  $A(t) \neq 0$, $B(t)$, $C(t) \neq 0$, $E(t)$, $F(t)$ and $Q(t)$ are arbitrary smooth functions of $t$.\\

We recall that an equivalence transformation of class (\ref{ed1}) is a nondegenerate transformation, which acts on independent and dependent variables given by

\begin{equation}
\label{equivalence} t= f(\tilde{t}, \tilde{x}, \tilde{u}), \quad x=g(\tilde{t}, \tilde{x}, \tilde{u}), \quad u=h(\tilde{t}, \tilde{x}, \tilde{u}),
\end{equation}

\noindent in a manner that it transforms equation (\ref{ed1}) in
\begin{eqnarray}
\label{eqtransformed}
\nonumber \begin{array}{c} \tilde{u}_{\tilde{t}}+ \tilde{A}(\tilde{t}) \tilde{u}_{\tilde{x}\tilde{x}\tilde{x}\tilde{x}\tilde{x}}+\tilde{B}(\tilde{t}) \tilde{u}_{\tilde{x}\tilde{x}\tilde{x}}+ \tilde{C}(\tilde{t}) \tilde{u} 
 \tilde{u}_{\tilde{x}\tilde{x}\tilde{x}}+ \tilde{E}(\tilde{t}) \tilde{u} \tilde{u}_{\tilde{x}}+ \tilde{F}(\tilde{t}) \tilde{u}_{\tilde{x}} \tilde{u}_{\tilde{x}\tilde{x}} + \tilde{Q}(\tilde{t}) \tilde{u}=0 ,
 \end{array}
\end{eqnarray}
i.e., it maps equation (\ref{ed1}) to another equation of the same class, where the transformed functions can be different from the original ones. The equivalence transformations have important applications. For instance, they can be used to determine equivalent formulations of a class of partial differential equations (PDEs) that could simplify the analysis of the class, especially if one expect to perform a classification problem. With respect to the concept of equivalence transformations, some results have been communicated in the recent literature \cite{cherniha,ndc,RGB:16,gandibra,garrido,gungor,ivanova,tortra2015,vaneeva}.\\

Nonlinear evolution equations (NLEEs) have been studied in a great number of works from the similarity reductions
(see e.g.  \cite{ademkhalique,anco,bokhari,freiremariano,liu,moitsheki,naz,rosa2016,sopho2011,wang}) which lead the NLEEs to ordinary differential equations (ODEs). However, it is not always evident how to integrate these ODEs. A different way is to consider another equation with a known solution, which is related with the equation under consideration through an invertible transformation. For a class of PDEs, this equivalence problem may be solved by using the differential invariants of its group of equivalence transformations.\\

A differential invariant of a class of PDEs  is a real valued function which is invariant under the equivalence group. The applications of differential invariants of the Lie groups of conti\-nuous transformations can be found in an extensive variety of problems arising in mathematical modeling, differential equations, differential geometry and nonlinear science. There exist different methods to approach the equivalence problem. Among them, the Lie infinitesimal method is one of the most frequently used. Based on Lie's preliminary results \cite{Lie}, Ovsyannikov \cite{ovsian} developed this approach for infinite Lie groups. Furthermore, this method was subsequently used for obtaining equivalence transformations. Afterwards, in \cite{olver}, Olver presented the equivariant method of moving frames which is used to generate differential invariants for finite dimensional Lie group actions. A broad theory on differential invariants of Lie groups along with algorithms of construction of differential invariants is given in \cite{olverinv,ovsian}.\\% (\ref{equivalence}). 

%For the general theory of differential invariants of Lie groups, see, e.g. Ovsiannikov LV. Group analysis of differential equations. New York: Academic Press; 1982..

Equivalence transformations play an important role in the theory of differential invariants. Indeed, the derivation of equivalence transformations of a class of PDEs %(\ref{ed1}) 
represents the initial step towards the determination of differential invariants. In \cite{ibrainv1997,ibrainv1}, Ibragimov proposed a method for calculation of differential invariants of classes of linear and nonlinear equations which have an infinite dimensional equivalence group. Following Ibragimov's method several authors have constructed differential invariants, see e.g. \cite{bagderina,bagderina14,mahomed,sophorita2008,tsaousi}. Furthermore, the knowledge of differential invariants can be used to find special equivalence transformations that linearize some classes of equations (see e.g.\cite{ibra2005,torrisi,rita}).\\

This work is organized as follows. In the next section,% we give a brief introduction to the continuous group of equivalence transformations of a PDE and 
we show the equivalence transformations of class (\ref{ed1}) and  of its subclass
\begin{equation}\label{ed1*}
\begin{array}{l}
 u_t+A(t) u_{xxxxx}+B(t) u_{xxx}+ C(t) u u_{xxx} +E(t) u u_x+ F(t) u_x u_{xx} =0. 
\end{array}
\end{equation}

 In Section \ref{diffinv}, following Ibragimov's method, we get differential invariants of zero and first order for class (\ref{ed1*}). As an application of the invariants obtained, in Section \ref{exact solutions}, we determine exact solutions for special subclasses of family (\ref{ed1}).

\section{Equivalence transformations}\label{equivalence transformations}
An equivalence transformation of class (\ref{ed1}) is a nondegenerate point transformation from $(t,x,u)$ to $(\tilde{t},\tilde{x},\tilde{u})$ that preserves the differential structure of the equation, that is, it maps every equation of class (\ref{ed1}) into another equation of the same class  but with different functions, $\tilde{A}(\tilde{t})$, $\tilde{B}(\tilde{t})$, $\tilde{C}(\tilde{t})$, $\tilde{E}(\tilde{t})$, $\tilde{F}(\tilde{t})$ and $\tilde{Q}(\tilde{t})$.
The set of equivalence transformations forms a group. %denoted by ${\cal G}$. 
Generators of the group of equivalence transformations of equation (\ref{ed1}) take the form

\begin{equation}\label{generator}
 Y  =   \tau \partial_t + \xi \partial_x+\eta \partial_u+ \omega_1 \partial_A  + \omega_2 \partial_B + \omega_3 \partial_C + \omega_4 \partial_E + \omega_5 \partial_F + \omega_6 \partial_Q .   \end{equation}	
								
In \cite{ndc} it was proved  that class (\ref{ed1}) admits an infinite dimensional equivalence group ${\cal G}$ whose generators are given by 
\begin{eqnarray}{}
%\label{eqgen1} \displaystyle Y_1 & = & x \partial_x+ \frac{u}{2} \partial_u + 5 A \partial_A + 3 B \partial_B + \frac{5 C}{2}  \partial_C  +  \frac{E}{2} \partial_E+ \frac{5 F}{2} \partial_F,\\
%\nonumber \\
\label{eqgen1} \displaystyle Y_1 & = & x \partial_x+  5 A \partial_A + 3 B \partial_B + 3C  \partial_C  +E \partial_E+ 3F \partial_F,\\
\nonumber \\
\label{eqgen2} \displaystyle Y_2 & = & \partial_x,\\
\nonumber \\
\label{eqgen3} Y_\alpha & = & \alpha \partial_t- \alpha_t A \partial_A- \alpha_t B \partial_B- \alpha_t C \partial_C - \alpha_t E \partial_E - \alpha_t F \partial_F- \alpha_t Q \partial_Q,\\
\nonumber \\
\label{eqgen4} Y_{\beta} & = &{\beta}  u \partial_u - {\beta}  C \partial_C- {\beta}  E \partial_E- {\beta}  F \partial_F- {\beta} _t \partial_Q,
\end{eqnarray}

\noindent where $\alpha=\alpha(t)$ and ${\beta} ={\beta} (t)$ are arbitrary smooth functions of $t$. From generators (\ref{eqgen1})-(\ref{eqgen4}) the finite form of the equivalence transformations was obtained 
\begin{equation}\label{transform}\begin{array}{c}
 \displaystyle \tilde{t}  =  \lambda(t), \quad \tilde{x}  =  (x+k_2)e^{k_1}, \quad \tilde{u}  =  e^{\frac{k_1}{2}+k_3 \mu(t) }u,
\end{array}\end{equation}
\noindent where $k_1$, $k_2$ and $k_3$ are arbitrary constants, $\lambda=\lambda(t)$ and $\mu=\mu(t)$ are arbitrary smooth functions with $\lambda_t \neq 0$.\\

\noindent We observe that the transformation
\begin{equation}\label{transiniz}\begin{array}{c}
 \displaystyle \tilde{t}  = t, \quad \tilde{x}  = x, \quad \tilde{u}  =  e^{\int{Q(t)dt}}u,
\end{array}\end{equation}
maps equations of the subclass
\begin{equation}
 \tilde{u}_{\tilde{t}}+\tilde{A}(\tilde{t}) \tilde{u}_{\tilde{x}\tilde{x}\tilde{x}\tilde{x}\tilde{x}}+
\tilde{B}(\tilde{t}) \tilde{u}_{\tilde{x}\tilde{x}\tilde{x}}+
 \tilde{C}(\tilde{t}) \tilde{u} \tilde{u}_{\tilde{x}\tilde{x}\tilde{x}}+\tilde{E}(\tilde{t}) \tilde{u} \tilde{u}_{\tilde{x}}+ 
\tilde{F}(\tilde{t}) \tilde{u}_{\tilde{x}} \tilde{u}_{\tilde{x}\tilde{x}} =0,
\end{equation}
in an equation of class \eqref{ed1}.

The subclass (\ref{ed1*}) admits an infinite dimensional equivalence group ${\cal G'}$ whose generators are given by 
\begin{eqnarray}{}
\label{eqgen1*} \displaystyle Y_1 & = & x \partial_x+  5 A \partial_A + 3 B \partial_B + 3C  \partial_C  +E \partial_E+ 3F \partial_F,\\
\nonumber \\
\label{eqgen2*} \displaystyle Y_2 & = & \partial_x,\\
\nonumber \\
\label{eqgen3*} Y_3 & = & u \partial_u - C \partial_C-  E \partial_E-  F \partial_F,\\
\nonumber \\
\label{eqgen4*} Y_\alpha & = & \alpha \partial_t- \alpha_t A \partial_A- \alpha_t B \partial_B- \alpha_t C \partial_C - \alpha_t E \partial_E - \alpha_t F \partial_F,
\end{eqnarray}
where $\alpha=\alpha(t)$ is an arbitrary smooth function of $t$.\\

 From generators (\ref{eqgen1*})-(\ref{eqgen3*}) the finite form of the equivalence transformations was obtained 
\begin{equation}\label{transform*}\begin{array}{c}
 \displaystyle
  \tilde{t}  =  \lambda(t), \quad \tilde{x}  =  (x+k_2)e^{k_1}, \quad \tilde{u}  =  e^{k_3 }u,
%  \\\nonumber \\
%\tilde{A}=\frac{e^{5k_1}}{\lambda_t}A ,\quad 
%\tilde{B}=\frac{e^{3k_1}}{\lambda_t}B ,\quad 
%\tilde{C}=\frac{e^{3k_1-k_3}}{\lambda_t}C ,\quad 
%\tilde{E}=\frac{e^{k_1-k_3}}{\lambda_t}E ,\quad 
%\tilde{F}=\frac{e^{3k_1-k_3}}{\lambda_t}F ,
\end{array}\end{equation}
\noindent where $k_1$, $k_2$ and $k_{3}$ are arbitrary constants, $\lambda=\lambda(t)$ is an arbitrary smooth function with $\lambda_t \neq 0$.

\section{Differential invariants of class (\ref{ed1*})}\label{diffinv}

Let us consider the continuous group ${\cal G'}$ of equivalence transformations of class (\ref{ed1*}) which is spanned by (\ref{eqgen1*})-(\ref{eqgen4*}). A differential invariant of zero order of class (\ref{ed1*}) is a real valued function $J$ of the form
\begin{equation}
\label{invariantzero} J(t,x,u,A,B,C,E,F),
\end{equation}
which is invariant under the equivalence group ${\cal G'}$. Similarly, we call differential invariant of order $s$ to an invariant function which has the form
\begin{eqnarray}
\nonumber \label{invariants} J(t,x,u,A,B,C,E,F, A_t, B_t, C_t, E_t, F_t,\ldots),
\end{eqnarray}
where $s$ denotes the maximal order derivative of the arbitrary functions $A$, $B$, $C$, $E$, and $F$.\\

In order to get the differential invariants of order $s$ of class (\ref{ed1*}) we need $Y^{(s)}$, that is the $s$th prolongation of the operator $Y$ given by (\ref{generator}).
%\begin{eqnarray}\label{prolongY}
%\nonumber \begin{array}{c}
%Y^{(1)}= Y+ \tilde{\omega}_1 \partial_{A_t}+\tilde{\omega}_2 \partial_{B_t}+\tilde{\omega}_3 \partial_{C_t} +\tilde{\omega}_4 \partial_{E_t}+\tilde{\omega}_5 \partial_{F_t}+\tilde{\omega}_6 \partial_{Q_t},
%\end{array}
%\end{eqnarray}
%
%\noindent with $\tilde{\omega}_i= \tilde{D}_t(\omega_i)-f^i_t \tilde{D}_t(\tau)$, where $f^i$, $i=1,\ldots,6$, represents each component $(A,B,C,E,F,Q)$. Here
%\begin{eqnarray}\begin{array}{c} \nonumber \tilde{D}_t = \partial_t + A_t \partial_A+ B_t \partial_B+...+ Q_t \partial_Q + A_{tt} \partial_{A_t}+ B_{tt} \partial_{B_t}+ \cdots+ Q_{tt} \partial_{Q_t}+ \cdots. 
%\end{array}\end{eqnarray}
For information on how the prolongations can be determined, the interested reader can refer, as example, to reference  \cite{sophorita2008} and references inside. 

\noindent Furthermore,  an equation 
\begin{eqnarray}\nonumber
H(t,x,u,A,B,C,E,F, A_t, B_t, C_t, E_t, F_t,\ldots)=0,
\end{eqnarray}
which satisfies the condition
\begin{eqnarray}
\nonumber \label{inveq} Y_{k}^{(s)}\left( H \right) \,{|}_{H=0}=0,
\end{eqnarray}
where $Y_k^{(s)}$ is the $s$th prolongation of the operator $Y_k$, is called invariant equation.\\

%
%The differential invariants of first order are determined by applying the first extensions of the operators (\ref{eqgen1})-(\ref{eqgen4}) on the function $J$ of the form
%\begin{eqnarray}
%\nonumber J(t,x,u,A,B,C,E,F,Q, A_t, B_t, C_t, E_t, F_t, Q_t).
%\end{eqnarray}
%Then the differential invariants of first order are determined by 
%the equations

%\begin{equation}
%\label{invcriterionfirstorder} Y_k^{(1)} (J)=0, \qquad k=1,2,\alpha,r.
%\end{equation}
%
%The invariance criterion (\ref{invcriterionfirstorder}) leads to a system of linear first order PDEs whose independent solutions provide the set of differential invariants of first order. 

\subsection{Differential invariants of zero order}
\medskip

Due to the fact that in equations (\ref{ed1*}) the arbitrary functions depend only on $t$, we search for invariant functions of the form
\begin{eqnarray}\label{inv0}
 J=J(t,A,B,C,E,F).
\end{eqnarray}

In order to get the differential invariants of zero order of class (\ref{ed1*}) we apply the operators (\ref{eqgen1*})-(\ref{eqgen4*}) on the function $J$, which is given by (\ref{inv0}), obtaining a system of linear first order PDEs for $J$:
\begin{eqnarray}
\nonumber \label{invcriterion} Y_k (J)=0, \qquad k=1,2,3,\alpha.
\end{eqnarray}
Solving this system
% we find the differential invariants of order zero. 
%\noindent Thus, invariant test $Y_2(J)=0$ is trivially satisfied. Using the operators (\ref{eqgen1}), (\ref{eqgen3}) and (\ref{eqgen4}), we obtain that $J$ must satisfy the following 
%equations
%\begin{equation}
%\label{opzeroorder} Y_1(J)=0, \quad Y_{\alpha}(J)=0, \quad Y_r(J)=0.
%\end{equation}
%
%\noindent Because $\alpha$, $\alpha_t$, $r$ and $r_t$ must be considered as independent functions, the invariant tests (\ref{opzeroorder}) lead to 5 first order linear PDEs for $J$. Taking into account these conditions, 
 we obtain two differential invariants of zero order
\begin{equation}\label{diffinv0}
J_1^{(0)} \displaystyle=\frac{A E}{B C}, \qquad J_2^{(0)}= \displaystyle \frac{F}{C},
\end{equation}
and the following invariant equations
\begin{equation}\label{inveqzero}
 A=0,\; B=0,\; C=0, E=0, F=0.
\end{equation}

\subsection{Differential invariants of first order}
\medskip

\noindent We are interested in obtaining differential invariants of first order
\begin{eqnarray}\label{inv1}
\nonumber J=J(t,A,B,C,E,F,A_t,B_t,C_t,E_t,F_t).
\end{eqnarray}
%The search of differential invariants of first order leads us to calculate the first prolongation of the operator $Y$ (\ref{generator}):
%\begin{eqnarray}\label{prolongY}
%\nonumber \begin{array}{c}
%Y^{(1)}= Y+ \tilde{\omega}_1 \partial_{A_t}+\tilde{\omega}_2 \partial_{B_t}+\tilde{\omega}_3 \partial_{C_t} +\tilde{\omega}_4 \partial_{E_t}+\tilde{\omega}_5 \partial_{F_t},
%\end{array}
%\end{eqnarray}
%
%\noindent with $\tilde{\omega}_i= \tilde{D}_t(\omega_i)-f^i_t \tilde{D}_t(\tau)$, where $f^i$, $i=1,\ldots,6$, represents each component $(A,B,C,E,F,Q)$. Here
%\begin{eqnarray}\begin{array}{c} \nonumber \tilde{D}_t = \partial_t + A_t \partial_A+ B_t \partial_B+...+ Q_t \partial_Q + A_{tt} \partial_{A_t}+ B_{tt} \partial_{B_t}+ \cdots+ Q_{tt} \partial_{Q_t}+ \cdots. 
%\end{array}\end{eqnarray}
%
%For further information on how the prolongations of higher order can be determined, one can refer, as example, to reference  \cite{sophorita2008}. 

\noindent The first prolongation of the operators (\ref{eqgen1*})-(\ref{eqgen4*}) yield the invariant criterion:

\begin{equation}
\label{opfirstorder} Y_1^{(1)}(J)=0,  \quad  Y_2^{(1)}(J)=0,  \quad  Y_3^{(1)}(J)=0, \quad Y_{\alpha}^{(1)}(J)=0.
\end{equation}
%
%\noindent Since $\alpha(t)$ and $r(t)$ are arbitrary functions, the invariant tests (\ref{opfirstorder}) lead us to a system of 7 first order linear PDEs for $J$. 
Solving this system we get % two differential invariants of zero order which are given by (\ref{diffinv0}) and 
the following differential invariants of first order
\begin{eqnarray}\label{diffinv1}
\begin{array}{rclcrcl}
\nonumber J_1^{(1)}&=& \displaystyle \frac{A\left(A_t B-A B_t\right)^2}{B^7}, &\quad &
J_2^{(1)}&=& \displaystyle \frac{A\left(A_t C- A C_t \right)^2}{B^5 \, C^2},\\ 
\nonumber J_3^{(1)}&=& \displaystyle \frac{A^3 \left(A_t E- A E_t \right)^2}{B^7 \, C^2}, &\quad &
J_4^{(1)}&=& \displaystyle \frac{A\left(A_t F- A F_t \right)^2}{B^5 \, C^2}.
\end{array}
\end{eqnarray}
Furthermore
\begin{equation}\label{inveq1}
\begin{array}{rcl}
 A_t B-A B_t &=& 0,    \\
 A_t C- A C_t  &=& 0,    \\ 
 A_t E- A E_t  &=& 0,    \\ 
 A_t F- A F_t &=& 0,   
\end{array}
\end{equation}
are invariant equations.\\

We observe that the equivalence group $\cal{G'}$ of class \eqref{ed1*} is a conditional equivalence group  (see e.g. \cite{van} and references inside) of the equivalence group $\cal{G}$ of class \eqref{ed1} under the condition $Q(t)=0$. It is evident that in this case  $\cal{G'}$ is a trivial conditional equivalence group because it is a subgroup of the equivalence group $\cal{G}$ obtained by setting $\beta(t)=1$ in the generator \eqref{eqgen4}. Then all differential invariants and invariant equations of class \eqref{ed1*} can be seen as invariants of class \eqref{ed1} with respect to the equivalence subgroup $\cal{G'}$.\\

Moreover it is possible to prove that two differential invariant of first order,  $J_3^{(1)}$ and $J_4^{(1)}$, can be obtained with the application of the invariant differentiation to the differential invariants of first order and further calculations show that the class of equations (\ref{ed1*}) admits four differential invariants of second order and four of third order.\\

 It seems  that the four differential invariants
\begin{equation}\label{basis}
J_1^{(0)} \displaystyle=\frac{A E}{B C}, \qquad J_2^{(0)}= \displaystyle \frac{F}{C},
\qquad
J_1^{(1)}= \displaystyle \frac{A\left(A_t B-A B_t\right)^2}{B^7},
\qquad J_2^{(1)}= \displaystyle \frac{A\left(A_t C- A C_t \right)^2}{B^5 \, C^2},
\end{equation}
constitute a basis for class (\ref{ed1*}). Consequently any higher order differential invariant can be obtained with the employment of the invariant differentiation. 
However this result, far of the spirit of this work,  needs to be proved and it will be a task for the near future.

\section{Exact solutions using differential invariants}
\label{exact solutions}

Two given PDEs are called equivalent if one can be transformed into the other by a change of variables. Once established the equivalence it is necessary to determine a transformation
that connects both partial differential equations. At this point it is possible to transform a known solution of a PDE to solutions of others which are equivalent to this one.
Any set of invariants can provide necessary conditions for deriving equivalent equations, since if two equations  are equivalent, they have the same differential invariants and they satisfy the same invariant equations.\\

Accordingly, taking into account the invariant equations (\ref{inveqzero}), we observe that any equation of class (\ref{ed1*}) with $B = 0$ cannot be connected by an equivalence transformation with an equation of class (\ref{ed1*}) with $B \neq 0$, because $B = 0$ is an invariant equation. The same occurs for $E = 0$ and $F = 0$. Moreover, we have already supposed that $A\neq 0$ and $C\neq 0$, then for this reason, we focus our attention on equations of class (\ref{ed1*}) with $ABCEF\neq 0$.\\

Here as an application of the results of the previous section, we look for the most general form of equations of class (\ref{ed1}) that can be mapped, under a point transformation of the form (\ref{transform}), to an equation of the same class of the form
\begin{equation}\label{eqcons}
\begin{array}{c}
u_t+ m_1 u_{xxxxx}+m_2 u_{xxx} + m_3 u u_{xxx} +m_4 u u_x+m_5 u_x u_{xx}=0,
\end{array}
\end{equation}
where $m_1$, $m_2$, $m_3$, $m_4$ and $m_5$ are arbitrary nonzero constants. 
First of all, we observe that equations of class (\ref{ed1}) can be mapped, under a point transformation of the form (\ref{transform}), to an equation of the special subclass (\ref{ed1*}).
Then we look for  equations of class (\ref{ed1*}) that can be mapped, under a point transformation of the form (\ref{transform*}), to an equation of the form \eqref{eqcons}.\\

\noindent For equations of subclass (\ref{eqcons}) is
\begin{eqnarray}
\begin{array}{c}
\nonumber A=m_1, \quad B=m_2, \quad  C=m_3, \quad E=m_4, \quad F=m_5.
\end{array}
\end{eqnarray}

\noindent For equations of subclass (\ref{eqcons})  the differential invariants  (\ref{basis}) assume the values
\begin{equation}\label{inveqcons}
J_1^{(0)}=\frac{m_1 m_4}{m_2 m_3}, \qquad J_2^{(0)}=\frac{m_5}{m_3},\qquad J_1^{(1)}=0,\qquad J_2^{(1)}=0.
\end{equation} 
Then in order to be transformed by an equivalence transformation of the form (\ref{transform*}), into equation (\ref{eqcons}), any equation of class (\ref{ed1*}) must assume the same values (\ref{inveqcons}) of the invariants $J_1^{(0)}$, $J_2^{(0)}$, $J_1^{(1)}$, and $J_2^{(1)}$, that is the following conditions must be satisfied
\begin{eqnarray}
\nonumber &&\displaystyle\frac{A E}{B C}=\frac{m_1 m_4}{m_2 m_3}, \\
\nonumber &&\displaystyle \frac{F}{C}=\frac{m_5}{m_3},\\
\nonumber && \displaystyle \frac{A\left(A B_t-A_t B\right)^2}{B^7}= 0, \\
\nonumber && \displaystyle \frac{A\left(A_t C- A C_t \right)^2}{B^5 \, C^2}  = 0. 
%\nonumber && A_t E- A E_t  = 0, \\ 
%\nonumber && A_t F- A F_t  = 0.
\end{eqnarray}
%Equations of subclass (\ref{eqcons}) satisfy the invariant equations (\ref{inveq1}), moreover for these equations the differential invariants of zero order (\ref{diffinv0}) assume the values
%\begin{equation}\label{inveqcons}
%J_1^{(0)}=\frac{m_1 m_4}{m_2 m_3}, \qquad J_2^{(0)}=\frac{m_5}{m_3}.
%\end{equation} 
%Then in order to be transformed by an equivalence transformation of the form (\ref{transform*}), into equation (\ref{eqcons}), any equation of class (\ref{ed1*}) must assume the same values (\ref{inveqcons}) of the invariants $J_1^{(0)}$ and  $J_2^{(0)}$, and must satisfy the same invariant equations (\ref{inveq1}), that is it must satisfy the following conditions
%\begin{eqnarray}
%\nonumber &&\displaystyle\frac{A E}{B C}=\frac{m_1 m_4}{m_2 m_3}, \\
%\nonumber &&\displaystyle \frac{F}{C}=\frac{m_5}{m_3},\\
%\nonumber && A B_t-A_t B = 0, \\
%\nonumber && A_t C- A C_t  = 0,   \\ 
%\nonumber && A_t E- A E_t  = 0, \\ 
%\nonumber && A_t F- A F_t  = 0.
%\end{eqnarray}
The more general form of equations of class (\ref{ed1*}), satisfying the previous conditions is
%\begin{equation}\label{ed2}
%\begin{array}{l}
%\displaystyle u_t + G(t)u_{xxxxx}+c_1 G(t)u_{xxx}+H(t) u u_{xxx} \\ \\ \qquad \displaystyle +\frac{c_1 m_1 m_4}{m_2 m_3}H(t) u u_{x}+ \frac{m_5}{m_3}H(t)u_x u_{xx} + \frac{G(t)H_t(t)-G_t(t)H(t)}{G(t)H(t)}u=0,
%\end{array}
%\end{equation}
\begin{equation}\label{ed2*}
\begin{array}{l}
\displaystyle u_t + G(t)u_{xxxxx}+c_1 G(t)u_{xxx}+\frac{c_2m_3}{m_5}G(t) u u_{xxx} 
%\\ \\ \qquad 
\displaystyle +\frac{c_1c_2 m_1 m_4}{m_2 m_5}G(t) u u_{x}+ c_2 G(t)u_x u_{xx} =0,
\end{array}
\end{equation}
where $c_1$, $c_2$ are arbitrary constants and $G(t) \neq 0$ is an arbitrary function. The equations (\ref{eqcons}) and (\ref{ed2*}) satisfy necessary conditions to be equivalent. It is simple to verify that the transformation 
\begin{equation}\label{transeq*}
\begin{array}{lll}
 \displaystyle
t  &=&  \displaystyle \left(\frac{c_1}{m_2}\right)^{\frac{5}{2}} m_1^{\frac{3}{2}} \int G(\tilde{t} ) d\tilde{t},\vspace*{0.05cm} \\
x  &=&  \displaystyle \sqrt{\frac{c_1 m_1}{m_2}}(\tilde{x}+k_2), \vspace*{0.05cm}\\
u  &=& \displaystyle  \frac{c_2m_2}{c_1 m_5}\tilde{u}, 
\end{array}
\end{equation}
 maps equations (\ref{eqcons}) to equations
\begin{equation}\label{transed2*} \begin{array}{l} \displaystyle
\tilde{u}_{\tilde{t}} + G(\tilde{t})\tilde{u}_{\tilde{x}\tilde{x}\tilde{x}\tilde{x}\tilde{x}}  +c_1 G(\tilde{t})\tilde{u}_{\tilde{x}\tilde{x}\tilde{x}}+\frac{c_2m_3}{m_5}G(\tilde{t}) \tilde{u} \tilde{u}_{\tilde{x}\tilde{x}\tilde{x}}% \\ \\ \displaystyle \qquad
 +\frac{c_1c_2 m_1 m_4}{m_2 m_5}G(\tilde{t}) \tilde{u} \tilde{u}_{\tilde{x}}+ c_2G(\tilde{t})\tilde{u}_{\tilde{x}} \tilde{u}_{\tilde{x}\tilde{x}}=0.
\end{array}
\end{equation}
Applying the transformation \eqref{transiniz}, that is by using the transformation
\begin{equation}\label{tras1s}\begin{array}{c}
 \displaystyle \tilde{t}  =  \bar{t}, \quad \tilde{x}  = \bar{x}, \quad \tilde{u}  =  e^{\int{Q(\bar{t})d\bar{t}}}\bar{u},
\end{array}\end{equation}
equation \eqref{transed2*} becomes
\begin{equation}\label{newtransf} \begin{array}{l} \displaystyle
\bar{u}_{\bar{t}} + G(\bar{t})\bar{u}_{\bar{x}\bar{x}\bar{x}\bar{x}\bar{x}}  +
c_1 G(\bar{t})\bar{u}_{\bar{x}\bar{x}\bar{x}}+\frac{c_2m_3}{m_5}G(\bar{t})e^{\int Q(\bar{t})\,d\bar{t} } \bar{u} \bar{u}_{\bar{x}\bar{x}\bar{x}} \\ \\ \displaystyle \qquad
 +\frac{c_1c_2 m_1 m_4}{m_2 m_5}G(\bar{t}) e^{\int Q(\bar{t})\,d\bar{t} }\bar{u} \bar{u}_{\bar{x}}+ c_2G(\bar{t})e^{\int Q(\bar{t})\,d\bar{t} }\bar{u}_{\bar{x}} \bar{u}_{\bar{x}\bar{x}}+Q(\bar{t}) \bar{u}=0.
\end{array}
\end{equation}
By introducing a function $H({\bar t}) $ such that  
\begin{equation}\label{H}
H(\bar{t})=\frac{c_2m_3}{m_5} G(\bar{t}) e^{\int Q(\bar{t})\,d\bar{t} },
\end{equation}
the transformed equation \eqref{newtransf} simplifies in
\begin{equation}\label{transed2} \begin{array}{l} \displaystyle
\bar{u}_{\bar{t}} + G(\bar{t})\bar{u}_{\bar{x}\bar{x}\bar{x}\bar{x}\bar{x}}  +c_1 G(\bar{t})\bar{u}_{\bar{x}\bar{x}\bar{x}}+H(\bar{t}) \bar{u} \bar{u}_{\bar{x}\bar{x}\bar{x}} \\ \\ \displaystyle \qquad +\frac{c_1 m_1 m_4}{m_2 m_3}H(\bar{t}) \bar{u} \bar{u}_{\bar{x}}+ \frac{m_5}{m_3}H(\bar{t})\bar{u}_{\bar{x}} \bar{u}_{\bar{x}\bar{x}} + \frac{G(\bar{t})H_{\bar{t}}(\bar{t})-G_{\bar{t}}(\bar{t})H(\bar{t})}{G(\bar{t})H(\bar{t})}\bar{u}=0.
\end{array}
\end{equation}
Summarizing, the composition of transformation \eqref{transeq*} with transformation \eqref{tras1s} with \eqref{H}
\begin{equation}\label{transeq} 
\begin{array}{lll}
\displaystyle 
t  &=&  \displaystyle \left(\frac{c_1}{m_2}\right)^{\frac{5}{2}} m_1^{\frac{3}{2}} \int G(\tilde{t} ) d\tilde{t},\vspace*{0.05cm} \\
x  &=&  \displaystyle \sqrt{\frac{c_1 m_1}{m_2}}(\tilde{x}+k_2), \vspace*{0.05cm}\\
u  &=& \displaystyle  \frac{m_2 H(\tilde{t})}{c_1 m_3 G(\tilde{t})}\tilde{u},
 \end{array}
 \end{equation}
 maps equations (\ref{eqcons}) to equations \eqref{transed2}.
Then transformation (\ref{transeq}) maps solutions of equations (\ref{eqcons}) into solutions of equations (\ref{transed2}).\\ 

We can look for travelling wave solutions for equations  (\ref{eqcons}), that is solutions of the form
\begin{eqnarray}
\nonumber u(t,x)=U(\sigma)\equiv U(x-at),
\end{eqnarray}
where $a$ is an arbitrary constant. Then the function $U$ must be solution of the ODE
\begin{eqnarray}\begin{array}{l}
\nonumber -aU'+ m_1U^{V}+m_2U''' 
+ m_3 UU'''  +m_4UU'+m_5U'U''=0.
\end{array}
\end{eqnarray}
We are able to obtain for this equation particular solutions, then in the following, by using transformation (\ref{transeq}), we show particular solutions for equations of the form (\ref{transed2}).

\begin{itemize}
\item[$\bullet$]
If $m_3\neq -m_5$, choosing\\ 

$ \begin{array}{l} a= \displaystyle \frac{m_4\left(m_1 m_4-m_2 m_5-m_2 m_3 \right)}{(m_3+m_5)^2}, \end{array}$
\smallskip

a solution of equation (\ref{eqcons}) is given by
\begin{eqnarray}\label{sol1} \begin{array}{lll}
\nonumber u(t,x)&=&\displaystyle c_2 \, \exp\left(\sqrt{\frac{-m_4}{m_3+m_5}} \left(x+  \frac{m_4\left(m_2 m_5+m_2 m_3-m_1 m_4  \right) t}{\left( m_3+m_5\right)^2}  \right)\right),
\end{array}
\end{eqnarray}
with $c_2$ an arbitrary constant. Thus, 
\begin{eqnarray}\label{transsol1}
\begin{array}{lll}
\nonumber \tilde{u}(\tilde{t},\tilde{x})&=& \displaystyle \frac{m_3 c_1 G(\tilde{t})}{m_2 H(\tilde{t})} \left( c_2 \exp\left(\sqrt{\frac{-m_4}{m_3+m_5}} \left( \left( \tilde{x} +k_2 \right)  \sqrt{\frac{c_1 m_1}{m_2}} \right.  \right. \right. \\ & & \displaystyle \left. \left. \left. + \frac{m_4\left(m_2 m_5+m_2 m_3-m_1 m_4  \right)}{\left( m_3+m_5\right)^2} \left( \frac{c_1}{m_2}\right)^{\frac{5}{2}} m_1^{\frac{3}{2}} \int G(\tilde{t}) d\tilde{t}
\right)\right) \right),
\end{array}
\end{eqnarray}
is a solution of equation (\ref{transed2}).\\

\item[$\bullet$] 
Consider $m_1$, $m_2$ and $m_3$ given by\\\\
$ \begin{array}{l} \displaystyle {m_1}={{-3\,a +2\,{ m_4}-2\,{ m_5}}\over{72}}, %\vspace*{0.1cm} 
\quad \displaystyle {m_2}={{-30\,a+17\,{ m_4}-8\,{m_5}}\over{36}}, %\vspace*{0.1cm} , \\ 
\quad\displaystyle{m_3}={{15\,a-10\,{m_4}+4\,{m_5}}\over{12}}, \end{array}$\\

%with $a$ arbitrary constant.
we obtain
\[
u(t,x)=\tanh ^2\left(x-a\,t\right),\]
 that it is a kink-type solution of  equation (\ref{eqcons}). Thus, 
$$\begin{array}{rcl} \tilde{u}(\tilde{t},\tilde{x})&=& \displaystyle \frac{3 c_1 \left(15 a-10 m_4+4 m_5\right) G(\tilde{t})}{\left(-30 a+17 m_4-8
   m_5\right) H(\tilde{t})}    \tanh ^2\left(9 \sqrt{2}  \, a \left(-3 a+2 m_4-2 m_5\right)^{\frac{3}{2}} \right.   \\ \\ 
   &&  \displaystyle \left(\frac{c_1}{-30 a+17 m_4-8 m_5}\right)^{\frac{5}{2}} \int G(\tilde{t}) \, d\tilde{t}  \left. -\left(\tilde{x}+k_2\right) \sqrt{\frac{c_1 \left(3 a-2 m_4+2
   m_5\right)}{60 a-34 m_4+16 m_5}}\right),\end{array} $$
   
\noindent is a solution of the following equation
$$\begin{array}{l}
\displaystyle \tilde{u}_{\tilde{t}} + G(\tilde{t})\tilde{u}_{\tilde{x}\tilde{x}\tilde{x}\tilde{x}\tilde{x}}+c_1 G(\tilde{t})\tilde{u}_{\tilde{x}\tilde{x}\tilde{x}} \\ \\ \qquad +H(\tilde{t}) \tilde{u} \tilde{u}_{\tilde{x}\tilde{x}\tilde{x}}  \displaystyle +\frac{6 c_1 m_4 \left(3 a-2 m_4+2 m_5\right) H(\tilde{t})}{\left(15 a-10 m_4+4 m_5\right) \left(30 a-17 m_4+8
   m_5\right)} \tilde{u} \tilde{u}_{\tilde{x}} \vspace*{0.1cm}\\ \qquad \quad \displaystyle + \frac{12 m_5 \, H(\tilde{t})}{15 a-10 m_4+4 m_5} \tilde{u}_{\tilde{x}} \tilde{u}_{\tilde{x}\tilde{x}} + \frac{G(\tilde{t})H_{\tilde{t}}(\tilde{t})-G_{\tilde{t}}(\tilde{t})H(\tilde{t})}{G(\tilde{t})H(\tilde{t})}\tilde{u}=0.
\end{array}$$

\item[$\bullet$] 
For $m_1$, $m_2$ and $m_3$ given by\\\\
$ \begin{array}{l} \displaystyle {m_1}={{-3\,c +{ m_4}+2\,{ m_5}}\over{72}}, %\vspace*{0.1cm} \\
\quad\displaystyle {m_2}={{15\,c-2\,{ m_4}-4\,{m_5}}\over{36}}, %\vspace*{0.1cm} , \\ 
\quad\displaystyle{m_3}={{-15\,c+5\,{m_4}+4\,{m_5}}\over{12}}, \end{array}$

we obtain a soliton solution of equation (\ref{eqcons})
\[u(t,x)={\rm sech}^2\left(x-a\,t\right).\]

Therefore,
$$\begin{array}{rcl} \tilde{u}(\tilde{t},\tilde{x})&=& \displaystyle \frac{3 c_1 \left(-15 a+5 m_4+4 m_5\right) G(\tilde{t})}{\left(15 a-2 m_4-4  m_5\right) H(\tilde{t})}   {\rm sech}^2\left(\left(\tilde{x}+  k_2\right)    \sqrt{\frac{c_1 \left(-3 a+m_4+2 m_5\right)}{30 a-4 m_4-8 m_5}}\right. \vspace*{0.1cm} \\ & & \displaystyle  -9
   \sqrt{2} \, a \left(-3 a+m_4+2 m_5\right)^{\frac{3}{2}} \vspace*{0.1cm}   \displaystyle \left. \left(\frac{c_1}{15 a-2 m_4-4 m_5}\right)^{\frac{5}{2}}
   \int G(\tilde{t}) \, d\tilde{t}\right), \end{array}$$
   
satisfies the following equation
$$\begin{array}{l}
\displaystyle \tilde{u}_{\tilde{t}} + G(\tilde{t})\tilde{u}_{\tilde{x}\tilde{x}\tilde{x}\tilde{x}\tilde{x}}+c_1 G(\tilde{t})\tilde{u}_{\tilde{x}\tilde{x}\tilde{x}}+H(\tilde{t}) \tilde{u} \tilde{u}_{\tilde{x}\tilde{x}\tilde{x}}  \\ \\ \displaystyle \qquad +\frac{6 c_1 m_4 \left(-3 a+m_4+2 m_5\right) H(\tilde{t})}{\left(15 a-2 m_4-4 m_5\right) \left(-15 a+5 m_4+4
   m_5\right)} \tilde{u} \tilde{u}_{\tilde{x}} \vspace*{0.1cm}\\ \qquad \quad \displaystyle+ \frac{12 m_5\, H(\tilde{t})}{-15 a+5 m_4+4 m_5} \tilde{u}_{\tilde{x}} \tilde{u}_{\tilde{x}\tilde{x}} + \frac{G(\tilde{t})H_{\tilde{t}}(\tilde{t})-G_{\tilde{t}}(\tilde{t})H(\tilde{t})}{G(\tilde{t})H(\tilde{t})}\tilde{u}=0.
\end{array}$$

\item[$\bullet$] 
Taking $m_1$, $m_2$, $m_3$ and $m_4$ given by
\\\\
$ \begin{array}{l} 
 {m_1}= \displaystyle-{{1}\over{64\,a^4}}, \quad {m_2}=\displaystyle-{{5}\over{16\,a^2}}, \quad \displaystyle 
{ m_3}=\displaystyle-2\,{m_5}, \quad {m_4}=\displaystyle-16\,a^2\,
{m_5},
\end{array}$\\

we get the compacton-type solution of equation (\ref{eqcons}),
\[
u(t,x)=\cos^4\left(a\left(x-t\right)\right).\]

Applying (\ref{transeq}) we obtain 

$$\begin{array}{rcl} \tilde{u}(\tilde{t},\tilde{x})&=& \displaystyle \frac{32 \,a^2\, c_1\, m_5\, G(\tilde{t})}{5 H(\tilde{t})} \cos ^4\left(\frac{1}{50}\sqrt{\frac{c1}{5}} \left(25  \left(\tilde{x}   +k_2\right)+4  c_1^2  \int G(\tilde{t}) \, d\tilde{t}\right)\right), \end{array}$$

\noindent which is a solution of the following equation

$$\begin{array}{l}
\displaystyle \tilde{u}_{\tilde{t}} + G(\tilde{t})\tilde{u}_{\tilde{x}\tilde{x}\tilde{x}\tilde{x}\tilde{x}}+c_1 G(\tilde{t})\tilde{u}_{\tilde{x}\tilde{x}\tilde{x}}+H(\tilde{t}) \tilde{u} \tilde{u}_{\tilde{x}\tilde{x}\tilde{x}}   \\ \\ \qquad \displaystyle+\frac{2 c_1 \, H(\tilde{t})}{5} \tilde{u} \tilde{u}_{\tilde{x}}- \frac{H(\tilde{t})}{2} \tilde{u}_{\tilde{x}} \tilde{u}_{\tilde{x}\tilde{x}} + \frac{G(\tilde{t})H_{\tilde{t}}(\tilde{t})-G_{\tilde{t}}(\tilde{t})H(\tilde{t})}{G(\tilde{t})H(\tilde{t})}\tilde{u}=0.
\end{array}$$  

\end{itemize}

\section{Conclusions}\label{conclusions}

In this paper, by using the continuous group of equivalence transformations of the family of fifth-order KdV equations (\ref{ed1}) and its special subclass \eqref{ed1*}, we determine the differential invariants, along with their invariant equations, of zero and first order. Based on the fact that equivalent equations have the same differential invariants and satisfy the same invariant equations,  we have obtained exact solutions of the wide class of equations (\ref{transed2}). First, we establish that an equation of class (\ref{ed1}) can be mapped by the equivalence transformations (\ref{transiniz}) to an equation of family (\ref{ed1*}).
After we prove that the equations \eqref{transed2*}  of subclass (\ref{ed1*}) are equivalent with respect to equivalence transformations (\ref{transeq*}) to  equations of family (\ref{ed1*}) which have a simpler form (\ref{eqcons}). Furthermore, we determine by composition the transformation which connects both equations (\ref{eqcons}) and the equations  (\ref{transed2}) of family \eqref{ed1}. Then, we construct exact solutions of the simpler equations (\ref{eqcons}). Some of these solutions have physical interpretations. Finally, by using transformation (\ref{transeq}) we find exact solutions of equations (\ref{transed2}).

\section*{Acknowledgments}

R. Tracin\`{a} acknowledges Universit\`{a} di Catania for financial support, FIR \textit{Charge Transport in Graphene and Low dimensional Structures}. R. de la Rosa and M.S. Bruz\'on acknowledge the financial support from Junta de Andaluc\'ia group FQM-201, Universidad de C\'adiz. M.S. Bruz\'on also acknowledges the support of DGICYT project MTM2009-11875 with the participation of FEDER. R. de la Rosa also expresses his sincere gratitude to the Plan Propio de Investigaci\'on de la
Universidad de C\'adiz.


\begin{thebibliography}{00}

%% \bibitem{label}
%% Text of bibliographic item

\bibitem{ademkhalique}
Adem, K.R., Khalique, C.M.: Exact solutions and conservation laws of Zakharov-Kuznetsov modified equal width equation with power law nonlinearity, Nonlinear Anal. RWA \textbf{13}, 1692--1702 (2012)

\bibitem{anco}
Anco, S.C., Ali, S., Wolf, T.: Symmetry analysis and exact solutions of semilinear heat flow in multi-dimensions, J. Math. Anal. Appl. \textbf{379}, 748--763 (2011) 


\bibitem{bagderina} Bagderina, Y.Y.: Invariants of a family of third-order ordinary differential equations, J. Phys. A: Math. Theor. \textbf{42}, 085204 (2009)


\bibitem{bagderina14} Bagderina, Y.Y., Tarkhanov, N.N.: Differential invariants of a class of Lagrangian systems with two degrees of freedom, J. Math. Anal. Appl. \textbf{410}, 733--749 (2014)

%\bibitem{bluman}
%Chaolu T., Bluman G.: An algorithmic method for showing existence of nontrivial non-classical symmetries of partial differential equations without solving determining equations, J. Math. Anal. Appl. \textbf{411}, 281--296 (2014)


\bibitem{bokhari}
Bokhari, A.H., Al Dweik, A.Y., Kara, A.H. et al.: A symmetry analysis of some classes of evolutionary nonlinear (2+1)-diffusion equations with variable diffusivity, Nonlinear Dyn \textbf{62}, 127--138 (2010)

%\bibitem{freire} Bozhkov, Y., Freire, I.L.: Symmetry analysis of the bidimensional Lane-Emden systems, J. Math. Anal. Appl. \textbf{388}, 1279--1284 (2012) 


%\bibitem{bruzon}
%Bruz\'on, M.S., Gandarias, M.L.: Symmetries for a family of Boussinesq equations with nonlinear dispersion, Commun Nonlinear Sci Numer Simulat \textbf{14}, 3250--3257 (2009)

%\bibitem{bgtt}
%Bruz\'on, M.S., Gandarias, M.L., Torrisi, M., Tracin\`{a}, R.: On some applications of transformation groups to a class of nonlinear dispersive equations, Nonlinear Anal. RWA \textbf{13}, 1139--1151 (2012) 


\bibitem{cherniha}

Cherniha, R., Serov, M., Rassokha, I.: Lie symmetries and form-preserving transformations of reaction-diffusion-convection equations, J. Math. Anal. Appl. \textbf{342}, 1363--1379 (2008) 


\bibitem{ndc}
de la Rosa, R., Gandarias, M.L., Bruz\'on, M.S.: Symmetries and conservation laws of a fifth-order KdV equation with time-dependent coefficients and linear damping, Nonlinear Dyn \textbf{84}, 135--141 (2016)

\bibitem{RGB:16}
de la Rosa, R., Gandarias, M.L., Bruz\'on, M.S.: Equivalence transformations and conservation laws for a generalized variable-coefficient Gardner equation, Commun Nonlinear Sci Numer Simulat \textbf{40}, 71--79 (2016)


\bibitem{freiremariano}
Freire, I.L., Torrisi, M.: Symmetry methods in mathematical modeling of \textit{Aedes aegypti} dispersal dynamics, Nonlinear Anal. RWA \textbf{14}, 1300--1307 (2013) 


\bibitem{gandibra}
Gandarias, M.L., Ibragimov, N.H.: Equivalence group of a fourth-order evolution equation unifying various non-linear models, Commun Nonlinear Sci Numer Simulat \textbf{13}, 259--268 (2008) 

\bibitem{garrido}  Garrido, T.M., Kasatkin, A.A., Bruz\'on, M.S., Gazizov, R.K.: Lie symmetries and equivalence transformations for the
Barenblatt-Gilman model, Journal of Computational and Applied Mathematics \textbf{318}, 253--258 (2017) 


 


%\bibitem{gungor}
%G\"{u}ng\"{o}r, F., Winternitz, P.: Equivalence classes and symmetries of the variable coefficient Kadomtsev-Petviashvili equation, Nonlinear Dyn \textbf{35}, 381--396 (2004)

\bibitem{gungor}
G\"{u}ng\"{o}r, F., \"{O}zemir, C.: Lie symmetries of a generalized Kuznetsov-Zabolotskaya-Khokhlov equation, J. Math. Anal. Appl. \textbf{423}, 623--638 (2015)

\bibitem{ibrainv1997}
Ibragimov NH.: Infinitesimal method in the theory of invariants of algebraic and differential equations, Not. S. Afr. Math. Soc. \textbf{29}, 61--70 (1997)

\bibitem{ibrainv1}
Ibragimov, N.H.: Elementary Lie group analysis and ordinary differential equations. John Wiley \& Sons, New York (1999)

\bibitem{ibra2005}

Ibragimov, N.H., Meleshko, S.V.: Linearization of third-order ordinary differential equations by point and contact transformations, J. Math. Anal. Appl. \textbf{308}, 266--289 (2005) 

\bibitem{ivanova}
Ivanova, N.M., Sophocleous, C.: On the group classification of variable-coefficient nonlinear diffusion-convection equations, Journal of Computational and Applied Mathematics \textbf{197}, 322--344 (2006) 


%\bibitem{ibrainv3}
%Ibragimov, N.H.: Laplace type invariants for parabolic equations, Nonlinear Dyn \textbf{28}, 125--133 (2002)

%\bibitem{ibrainv2}
%Ibragimov, N.H.: Invariants of a remarkable family of nonlinear equations, Nonlinear Dyn \textbf{30}, 155--166 (2002)

%\bibitem{khalique}
%Johnpillai, A.G., Khalique, C.M.: Lie group classification and invariant solutions of mKdV equation
%with time-dependent coefficients, Commun Nonlinear Sci Numer Simulat \textbf{16}, 1207--1215 (2011) 

\bibitem{Lie}

Lie, S.: \"{U}ber Differentialinvarianten, Math. Ann. \textbf{24}, 52--89 (1884)

\bibitem{liu}
Liu, H., Li, J.: Symmetry reductions, dynamical behavior and exact explicit solutions to the Gordon types of equations, Journal of Computational and Applied Mathematics \textbf{257}, 144--156 (2014) 

\bibitem{mahomed}
Mahomed, F.M., Johnpillai, A.G., Aslam, A.: Symmetry classification and joint invariants for the scalar linear
$(1 + 1)$ elliptic equation, Commun Nonlinear Sci Numer Simulat \textbf{25}, 84--93 (2015) 

\bibitem{moitsheki}
Moitsheki, R.J., Hayat, T., Malik, M.Y.: Some exact solutions of the fin problem with a power law
temperature-dependent thermal conductivity, Nonlinear Anal. RWA \textbf{11}, 3287--3294 (2010) 

\bibitem{naz}

Naz, R.: Group invariant solution for a free jet on a hemi-spherical shell, Applied Mathematics and Computation \textbf{215}, 3265--3270 (2010) 

\bibitem{olverinv}
Olver, P.: Equivalence, invariants and symmetry. Cambridge University Press, Cambridge (1995)

\bibitem{olver}
Olver, P.: Generating differential invariants, J. Math. Anal. Appl. \textbf{333}, 450--471 (2007) 


\bibitem{ovsian}
Ovsyannikov, L.V.: Group analysis of differential equations. Academic Press, New York (1982)

\bibitem{rosa2016}
Rosa, M., Camacho, J.C., Bruz\'on, M.S., Gandarias, M.L.: Classical and potential symmetries for a generalized Fisher equation, Journal of Computational and Applied Mathematics \textbf{318}, 181--188 (2017) 


\bibitem{sophorita2008}
Sophocleous, C., Tracin\`{a}, R.: Differential invariants for quasi-linear and semi-linear wave-type equations, Applied Mathematics and Computation \textbf{202}, 216--228 (2008) 

\bibitem{sopho2011}
Sophocleous, C., O'Hara, J.G., Leach, P.G.L.: Symmetry analysis of a model of stochastic volatility with
time-dependent parameters, Journal of Computational and Applied Mathematics \textbf{235}, 4158--4164 (2011) 




\bibitem{torrisi}
Torrisi, M., Tracin\`{a}, R., Valenti, A.: On the linearization of semilinear wave equations, Nonlinear Dyn \textbf{36}, 97--106 (2004)

%\bibitem{tortra2011}
%Torrisi, M., Tracin\`{a}, R.: Exact solutions of a reaction-diffusion system for \textit{Proteus mirabilis}
%bacterial colonies, Nonlinear Anal. RWA \textbf{12}, 1865--1874 (2011) 

\bibitem{tortra2015}
Torrisi, M., Tracin\`{a}, R.: An Application of Equivalence Transformations to Reaction Diffusion Equations, Symmetry \textbf{7}, 1929--1944 (2015)



\bibitem{rita}
Tracin\`{a}, R.: Invariants of a family of nonlinear wave equations, Commun Nonlinear Sci Numer Simulat \textbf{9}, 127--133 (2004)

\bibitem{tsaousi}

Tsaousi, C., Tracin\`{a}, R., Sophocleous, C.: Differential invariants for third-order evolution equations, Commun Nonlinear Sci Numer Simulat \textbf{20}, 352--359 (2015) 



%\bibitem{vanejohn}
%Vaneeva, O., Johnpillai, A.G., Popovych, R.O., Sophocleous, C.: Enhanced group analysis and conservation laws
%of variable coefficient reaction-diffusion equations with power nonlinearities, J. Math. Anal. Appl. \textbf{330}, 1363--1386 (2007) 

%\bibitem{vaneeva} Vaneeva, O., Kuriksha, O., Sophocleous, C.: Enhanced group classification of Gardner equations with time-dependent coefficients, Commun Nonlinear Sci Numer Simulat \textbf{22}, 1243--1251 (2015) 

\bibitem{van} 
Vaneeva O., Johnpillai A.G.,Popovych  R.O., Sophocleous C.: Enhanced group analysis and conservation laws of variable coefficient reaction-diffusion equations with power nonlinearities, J. Math. Anal. Appl. \textbf{330},  1363--1386 (2007)

\bibitem{vaneeva}
Vaneeva O., Popovych R.O., Sophocleous C.: Extended group analysis of variable coefficient reaction-diffusion equations with exponential nonlinearities, J. Math. Anal. Appl. \textbf{396}, 225--242 (2012) 

\bibitem{wang}  Wang, G., Kara, A.H., Fakhar, K.: Nonlocal symmetry analysis and conservation laws to an third-order Burgers equation, Nonlinear Dyn \textbf{83}, 2281--2292 (2016)





\end{thebibliography}
\end{document}